\documentclass{amsart}

\usepackage{amsfonts}
\usepackage{amssymb}
\usepackage{amsmath,amsthm}
\usepackage{stmaryrd}

\usepackage{xcolor} 

\usepackage[hidelinks]{hyperref} 
\hypersetup{
    colorlinks=true,
    linkcolor=black,
    citecolor=blue,
    filecolor=black,
   urlcolor=black}

\newtheorem{theorem}{Theorem}[section]
\newtheorem{definition}{Definition}[section]
\newtheorem{lemma}[theorem]{Lemma}

\newtheorem*{maintheorem*}{Theorem}

\newcommand{\fullgp}[1]{\llbracket #1 \rrbracket}

\newcommand\restr[2]{{
  \left.\kern-\nulldelimiterspace 
  #1 
  \littletaller 
  \right|_{#2} 
  }}
\newcommand{\littletaller}{\mathchoice{\vphantom{\big|}}{}{}{}}
\begin{document}

\title{Finitely Presented Covers of Topological Full Groups }
\author{Mustafa G\"{o}khan Benl\.{I}}

\address{\parbox{\linewidth}{
  Mustafa G\"{o}khan Benl\.{I} \\
  Department of Mathematics, Middle East Technical University\\
  Ankara, Turkey
}}
\email{benli@metu.edu.tr}

\author{Anıl \"{O}zdemir}
\address{\parbox{\linewidth}{
 Anıl \"{O}zdemir \\
  Department of Mathematics, TED University\\
  Ankara, Turkey}
}
\email{anil.ozdemir@tedu.edu.tr}


\begin{abstract} 
Commutator subgroups of infinite  minimal subshifts provide the first examples of infinite, finitely generated, simple amenable groups.
We  prove that if $G$ is a finitely presented group which maps onto 
the commutator subgroup of the topological full group of an infinite minimal subshift, then $G$ contains a non-abelian free subgroup. This in particular implies that these groups are not finitely presented. We also show that the topological full group of an irreducible subshift of finite type contains non-abelian free subgroups.

\end{abstract}

\maketitle
\tableofcontents
\section{Introduction}
Topological full groups stand at a fascinating intersection of topological dynamics, geometric group theory and operator algebras. They were studied by Giordano, Putnam and Skau \cite{GPS1999}  to classify Cantor minimal systems up to flip conjugacy or orbit equivalence.
Roughly, the topological full group $\fullgp{T}\le \mathrm{Homeo}(X)$ of a Cantor system  
$(X,T)$ consist of homeomorphisms that locally act like a power of $T$. Remarkably, the algebraic structure of the topological full group often completely encodes the underlying dynamics into a group-theoretic invariant. 

Beyond their role in classifying Cantor minimal system, topological full groups have emerged as a fertile ground for constructing exotic infinite groups, yielding the first examples of infinite, finitely generated, simple, amenable groups whose existence was a long-standing conjecture in modern group theory (See \cite{Matui2006}, \cite{JM2013}).

Algebraic properties of topological full groups were first investigated by Matui (\cite{Matui2006,Matui2013,Matui2015}). He proved that the commutator subgroup of these groups are simple and
showed that the commutator subgroup is finitely generated exactly when the system is a minimal subshift. He also showed that these groups are never finitely presented. Later Grigorchuk and Medynets investigated these groups further and showed that these groups are locally embeddable into finite groups and also  provided an infinite presentation of these groups (\cite{GM2014,GM2018}). In their breakthrough paper (\cite{JM2013}) Jushenko and Monod showed that these groups are amenable, thus constructing the first examples of infinite, finitely generated, simple,  groups.

It was shown by Matte Bon (\cite{MB2014}) that groups of intermediate growth constructed by
Grigorchuk (\cite{Grigorchuk1984}) can be embedded into commutator subgroups topological full groups. Also, the full group of a substitutional subshift containing the \textit{first Grigorchuk group} of intermediate growth was investigated in detail by Vorobets in \cite{Vorobets2020}.

The motivating question of this article is the following: Given a minimal subshift $(X,T)$, which finitely presented groups
can cover the finitely generated group $\fullgp{T}'$? (We say a group $G$\textit{ covers} a group $H$ if $G$ surjects onto $H$). This question is meaningful from two viewpoints: First, these groups are not  finitely presented. Secondly,  there are finitely generated amenable groups (or groups of intermediate growth) for which every finitely presented cover contains non-abelian free groups (see \cite{BS1980,BGdlH2013}). It is an open question if there is a finitely generated group of intermediate growth which has a finitely presented amenable cover (see \cite{Juschenko2022}). Until very recently, the existence of finitely presented, infinite, simple, amenable groups was open. One of the many results shared by Open AI recently claims 
to construct the first examples such groups (\cite{openai2026sharing}).

\vspace{0.2cm}
The main theorem of this article (Theorem \ref{mainthm})  is the following:

\begin{maintheorem*}
Let $(X,T)$ be an infinite minimal subshift. If $G$ is a finitely presented group
which maps onto $\fullgp{T}'$, then $G$ contains a non-abelian free subgroup. In 
particular, $\fullgp{T}'$ is not finitely presented.
\end{maintheorem*}

This shows that commutator subgroups of topological full groups of minimal subshfits are more then being infinitely presented: They cannot be even covered with finitely presented groups which do not have non-abelian free subgroups. 

The idea of the main theorem is to  approximate a minimal subshifts by \textit{subshifts of finite type} which are subshifts of very different nature. This approximation translates 
to the plane of full groups as an approximation of an amenable group by non-amenable ones.
During this process we prove the following theorem (Theorem \ref{thm:irreducible-sft-free-subgroup}) which we think is of independent interest:

\begin{maintheorem*}
	
	Let $(X,T)$ be an infinite irreducible subshift of finite type. Then
	$\fullgp{T}'$ contains a non-abelian free subgroup.
\end{maintheorem*}

\textbf{Acknowledgements:} The authors wish to thank R.Grigorchuk,  H. Matui and K. Medynets for helpful private discussions regarding the results. No AI was used in the writing or proving the results.

\section{Preliminaries}

In this section we gather required preliminaries which will be used in the sequel of the article.

\subsection{Cantor Systems and their topological full groups} \mbox{}

\vspace{0.2cm}

Let $\mathcal{C}$ denote the Cantor set. Recall that $\mathcal{C}$ is (up to homeomorphism) the unique zero-dimensional, compact metric space which has no isolated points. Thus $\mathcal{C}$ has many models such as $2^\mathbb{N},2^{\mathbb{Z}}$ or $A^\mathbb{Z}$ where $A$ is a finite set with at least two elements. Also recall that any compact, metric space is a continuous image of $\mathcal{C}$.  By a \textit{Cantor system} we mean a tuple $(X,T)$ where $X$ is (homeomorphic to) the Cantor set and  $T:X\to X$ is a homeomorphism. The main object of study of \textit{Cantor dynamics} is such tuples,  viewed as topological  dynamical systems. Two such systems $(X,T)$ and $(Y,S)$ are \textit{conjugate} if there is a homeomorphism $\pi: X\to Y$
such that $S\circ \pi= \pi \circ T$.

Here is the main object of study of this article:

\begin{definition}
Let $(X,T)$ be a Cantor system. The topological full group  of this system consists of
all homeomorphisms $g:X\to X$ such that there is a clopen partition $C_1,\ldots,C_k$ of $X$ and integers
$n_1,\ldots,n_k\in \mathbb{Z}$ such that $\restr{g}{C_i}=\restr{T^{n_i}}{C_i}$ for all $i=1,\ldots,k$.
The topological full group of the system will be denoted by $\fullgp{T}$.
\end{definition} 

Equivalently, a homemorphism $g\in Homeo(X)$ belongs to $\fullgp{T}$ if and only if there is a continuous function
$n_g : X\to \mathbb{Z}$  such that for all $x\in X$, $g(x)=T^{n_g(x)}(x)$ (here $\mathbb{Z}$ has the discrete topology). In words, $\fullgp{T}$ consists of all homeomorphisms of which locally are powers of $T$.

If $(X,T)$ is a Cantor system and $x\in X$, the $T$-orbit of $x$ is $\{T^n(x) \mid n\in \mathbb{Z}\}$. The Cantor system is \textit{minimal} is every orbit is dense in $X$. This is equivalent to $X$ having no non-trivial $T$-invariant closed subsets. Minimality is an irreducibility condition on Cantor systems.

From the point of view of Cantor dynamics, the importance of topological full groups is that they are invariants up to \textit{flip conjugacy}: (Two Cantor systems $(X,T),(Y,S)$ are \textit{flip conjugate} if 
$T$ is conjugate to $S$ or $S^{-1}$.)

\begin{theorem}\cite[Corollary 4.6]{GPS1999}
 Let $(X,T),(Y,S)$ be two minimal Cantor systems. Then $\fullgp{T}$ and $\fullgp{S}$ are isomorphic if and only if the systems are flip conjugate.
\end{theorem}

A rich source of Cantor minimal systems come from symbolic dynamics. Let $A$ be a finite alphabet with at least two elements, so that $A^\mathbb{Z}$ is a Cantor set. The left shift $T:A^\mathbb{Z}\to A^\mathbb{Z}$, defined by $T(x)_{n}=x_{n+1}$, is a homeomorphism. A \textit{subshift} is a $T$-invariant closed subset $X\subseteq A^\mathbb{Z}$, thus $(X,T)$ is a  dynamical system, which is a Cantor system if it is infinite and minimal. There are uncountably many non-conjugate minimal subshifts (see \cite{MorseHedlund1940,Kurka2003}) and hence there are uncountable many Cantor minimal systems whose topological full groups are non-isomorphic.

Recall that a group $G$ is LEF (\textit{locally embeddable into finite groups}) if for every finite subset $F\subseteq G$ there is a finite group $H$ and an injective map $\phi:F\to H$ such that
whenever $x,y,xy\in F$, then $\phi(xy)=\phi(x)\phi(y)$. Note that residually finite groups are LEF and
LEF groups are sofic.

Regarding  group theoretic properties of topological full groups  we have the following list:

\begin{theorem} \label{ref} Let $(X,T)$ be a Cantor minimal system. Then

\begin{enumerate} 
\item The commutator subgroup $\fullgp{T}'$ of the full 
group $\fullgp{T}$ is a simple group.
\item  $\fullgp{T}'$ is finitely generated if and only if $(X,T)$ is conjugate to a minimal subshift.
\item  $\fullgp{T}$ is LEF.
\item  $\fullgp{T}'$ is the product of two infinite locally finite subgroups.
\item  $\fullgp{T}$ is amenable.
\item If $(X,T)$ is a minimal subshift, then the finitely generated group $\fullgp{T}'$ is not finitely presented.
\item If $(X,T)$ is a minimal subshift, then the finitely generated group $\fullgp{T}'$ has exponential growth.
\item For every $\omega\in \{0,1,2\}^\mathbb{N}$, there is a minimal subshift $(X_\omega,T_\omega)$ such that
the Grigorchuk group $G_\omega$ embeds into $\fullgp{T_\omega}'$.
\end{enumerate}

\end{theorem}

\begin{proof} 
(1), (2) and (6) are proven in  \cite{Matui2006} (Theorems 4.9, 5.4 and 5.7 respectively). Items (3) and (4) are proven in \cite{GM2014} (Theorems 5.1 and 2.10). Item (5) is Theorem A in  \cite{JM2013}. Item (6) was proven in \cite[Theorem 2.4]{Matui2013}. Also, it can be observed that finitely presented LEF groups are residually finite and hence items (1) and (3) imply item (6). This argument is given in \cite[Corollary 2.9]{GM2014}.   Grigorchuk groups were consturcuted in \cite{Grigorchuk1984} and contain the first examples of groups of intermediate growth. Item (7) is \cite[Theorem 2.4]{Matui2013}. Item (8)  is proven in \cite[Theorem 1.1]{MB2014}.
\end{proof}

Combining the above results gives uncountably many examples of infinite, finitely generated, simple, amenable groups. More precisely, the commutator subgroup of the topological full group of an infinite  minimal subshift is infinite, finitely generated, simple and amenable. These are the first knowns examples of such groups and thus have attracted attention of group theorists. Such groups also may contain groups of intermediate growth.
Note that these groups are not finitely presented by item (6) of the previous theorem.

\subsection{Symbolic Dynamics  } \label{symbdynand}\mbox{}

\vspace{0.2cm}

 Let $A$ be a finite alphabet. An arbitrary element $\in A^\mathbb{Z}$ will be written as a bi-infinite sequence  $x= \ldots x_{-2}x_{-1}\bullet x_{0}x_{1}\ldots$ where $\bullet$ signifies the position between coordinates $-1$ and $0$ of $x$. With this notation the left shift map $T:A^\mathbb{Z}\to A^\mathbb{Z}$ becomes $$T(\ldots x_{-2}x_{-1}\bullet x_{0}x_{1}\ldots)=\ldots x_{-2}x_{-1}x_{0}\bullet x_{1}\ldots$$
For $n\ge 1$, let $A^n$ denote the set of finite words over $A$ of length $n$ and let $A^\ast=\bigcup_{n\ge 1} A^n$ be the set of all finite words over $A$. If $w\in A^n$, the length of $w$ is $|w|=n$.   Also, for $x\in A^\mathbb{Z}$ and $k\in\mathbb{Z},\ell\ge 0$,  $x_{[k,k+\ell]}=x_kx_{k+1}\cdots x_{k+\ell-1}\in A^\ast$.

Given two  words $u=u_nu_{n-1}\cdots u_{1}$ and $w=w_0\cdots w_{m-1} \in A^\ast$, sets of the form 
$$[u\bullet w]=\{x\in A^\mathbb{Z}\mid x_{-i}=u_{i},1\le i \le n , x_j=w_j, 0\le j\le m-1\}\subseteq A^\mathbb{Z}$$
are called \textit{cylinder sets} and they form a clopen  basis for the topology on $A^\mathbb{Z}$.
If $x\in [u\bullet w]$ we will write $x=\ldots u \bullet w\ldots$.

 If $X\subseteq A^\mathbb{Z}$ is a subshift and $u,v\in A^\ast$, we again write $[u\bullet w]$ for the corresponding cylinder set in $X$ (which may be empty). For each  $n\ge 1$, let $L_n(X)\subseteq A^n$ be the set of words of length $n$ that appear in sequences of $X$ and let us call $L(X)=\bigcup_{n\ge 1}L_n(X)\subseteq A^\ast$  the \textit{language} of $X$.
 
 The following is well known:

\begin{lemma}\label{norepeat}\cite[Proposition 2.3 ]{GM2018} Let $(X,T)$ be a minimal subshift over an alphabet $A$. Then there is a minimal subshift $(Y,S)$ (over  a larger alphabet)  such that $(X,T)$ and $(Y,S)$ are conjugate
and $(Y,S)$ satisfies the following:

Any  word $w\in L_5(Y)$ has no repeated letters, i.e., if $w=w_0w_1w_2w_3w_4$, all $w_0,w_1,w_2,w_3,w_4$ 
are distinct. 

\end{lemma}

It is clear that conjugate subshifts have isomorphic full groups. Thus, in all what follows, given a minimal subshift we will always pass to a conjugate subshift which satisfies the conclusion Lemma \ref{norepeat}.

For an arbitrary Cantor minimal system $(X,T)$, let $U\subseteq X$ be clopen such that $U,T(U),T^2(U)$ are mutually disjoint. Define the homeomorphism $\sigma_U$ as follows:

$$\sigma_U(x)=\left\{\begin{array}{cc}
T(x) & \text{if} \;  x\in U \cup T(U) \\
T^{-2}(x) & \text{if} \;  x \in T^2(U) \\
x & \text{otherwise} 
\end{array}\right. $$
In other words, $\sigma_U$ permutes $U,T(U),T^2(U)$ like a $3$-cycle and  is clearly in $\fullgp{T}$.
Let us call elements of the form $\sigma_U$ a \textit{local 3 cycle} at $U$.

Now suppose $(X,T)$ be a minimal subshift  that satisfies the conclusion of Lemma  \ref{norepeat}. Given $u=abc\in L_3(X)$, note that $[a\bullet bc],T[a\bullet bc],T^2[a\bullet bc]$ are pairwise disjoint and hence $\sigma_{[a\bullet bc]}$ is well  defined. 

\begin{theorem}\label{mat_fin_gen}\cite[Theorem 5.4]{Matui2006},\cite[Theorem 3.8]{GM2018}
Let $(X,T)$ be a minimal subshift  that satisfies the conclusion of Lemma  \ref{norepeat}. Then $\fullgp{T}'$
is generated by $\sigma_{[a\bullet bc]}, abc \in L_3(X)$.
\end{theorem}

\vspace{0.5cm}

In \cite{GM2018}, Grigorchuk and Medynets provided an infinite presentation for $\fullgp{T}'$ where $(X,T)$ is a minimal subshift. We will now describe this presentation which is an important ingredient in the proof of our main theorem.

Let $(X,T)$ be a minimal subshift over an alphabet $A$ that satisfies the conclusion of Lemma  \ref{norepeat}. For $w\in L(X)$ and $i\in \mathbb{Z}$, let $(w,i)=T^i[\bullet w]$. Note that $\fullgp{T}'$ is generated by $\{\sigma_{(w,1)}\mid w\in L_3(X)\}$ by Theorem \ref{mat_fin_gen}. Let $S=\{ x_{(w,1)} \mid w\in L_3(X) \}$, these will be the generators of the presentation. Let also $F_S$ denote the free group on $S$. If $w\in A^3 \setminus L_3(X)$, let $x_{(w,1)}$ be the empty word in $F_S$.  

$$\text{For}\; a\in A,\; \text{let}\;\displaystyle x_{(a,1)}=x_{[a\bullet ]}=\prod_{b,c\in A}x_{(abc,1)}$$

$$\displaystyle \text{For}\; a,b\in A,\; \text{let}\; x_{(ab,1)}=\prod_{c\in A}x_{(abc,1)}
$$

$$\displaystyle \text{For}\; b,c\in A,\; \text{let}\; x_{[\bullet bc]}=\prod_{a\in A}x_{(abc,1)}
$$

(Note that similar identities  hold in $\fullgp{T}'$ when $x$ is replaced by $\sigma$.)

Thus, for every $\displaystyle w\in \bigcup_{i=1}^3 L_i(X)$ we have defined the words $x_{(w,1)}\in F_S$.
For any two elements $x,y$ in a group, let us use the notation $x \ast y= [x^{-1},y]=xy^{-1}x^{-1}y$.

Given $w=w_0w_1\cdots w_{n-1}\in L_n(X)$ with $n\ge 4$ define
$$ x_{(w,1)}=x_{[\bullet w_0w_1]}
\ast (x_{[\bullet w_1w_2]} 
\ast 
\cdots 
\ast
(x_{[\bullet w_{n-4}w_{n-3}]}\ast
x_{(w_{n-3}w_{n-2}w_{n-1},1)})\cdots ) $$

For example, if $n=4$ and  $w=abcd$ then
 $$x_{(abcd,1)}=x_{[\bullet ab]}\ast x_{(bcd,1)}= \displaystyle \left( \prod_{e\in A}x_{(eab,1)} \right) \ast x_{(bcd,1)}$$


(By \cite[Proposition 3.5]{GM2018} similar identities hold in $\fullgp{T}'$ when $x$ is replaced by $\sigma$.) 

So far, for each $w\in L(X)$ we have defined a word $x_{(w,1)}\in F_S$.
Now, given $w=w_0w_1\cdots w_{n-1}\in L_n(X)$ with $n\ge 4$ and $2\le k \le n-1$ define
$$x_{(w,k)}=  (\cdots((x_{(w,1)}\ast x_{[w_2\bullet]})\ast x_{[w_3\bullet]})\ast \cdots \ast x_{[w_k\bullet]}) $$
(Again, by \cite[Proposition 3.5]{GM2018} similar identities hold in $\fullgp{T}'$ when $x$ is replaced by $\sigma$.) Thus, for each $w\in L(X)$ and $1\le k \le |w|-1$, we have defined a word $x_{(w,k)}\in F_S$.

\begin{theorem}\label{gmpres} \cite[Theorem 4.6]{GM2018} Let $(X,T)$ be a minimal subshift satisfying Lemma  \ref{norepeat}. Then $\fullgp{T}'$ has the following presentation:

\vspace{0.2cm}
Generators: 
$S=\{ x_{(w,1)} \mid w\in L_3(X) \}$.

Relations:  
\begin{enumerate}
\item For any $w\in L(X)$ and $1\le i \le |w|-1$,  $\left( x_{(w,i)} \right)^3=1$.
\item For any  $w\in L(X)$ and $1\le i < |w|-1$,  $\left( x_{(w,i)}\cdot x_{(w,i+1)} \right)^2=1$.
\item For any $w\in L(X)$ and $1\le i < |w|-2$,  $ x_{(w,i+1)}= x_{(w,i)}*x_{(w,i+2)} $.
\item For any $w,v\in L(X)$ and $ 1\le i \le |w|-1,1\le j \le |v|-1$ 
$$[x_{(w,i)},x_{(v,j)}]=1$$ 
whenever  $(w,i),(w,i+1),w(i+2),(v,j),(v,j+1),(v,j+2)$ are disjoint.

\item  For any $w\in L(X)$ and $1\le i \le |w|-1$ and any partition $C$ of $(w,i)$
into cylinder sets of the form $(u,k)$ where $u\in  L(X)$ and $1\le k \le |u|-1$,
$$x_{(w,i)}= \prod_{(u,k)\in C} x_{(u,k)}. $$
\end{enumerate}
\end{theorem}
Relations (1)-(4) are reminiscent of the relations satisfied by $3$-cycles in alternating groups (e.g., relations of type (2) correspond to relations of the form $((i,i+1,i+2)(i+1,i+2,i+3))^2=1$) . Relations of type (5) correspond to the general fact that if a clopen set is $U$ is the disjoint  of $U_1,U_2$ then $\sigma_U=\sigma_{U_1}\sigma_{U_2}$.

\vspace{0.5cm}

Another important class of subshifts are \textit{shifts of finite type}:

\begin{definition}
Let $A$ be a finite alphabet and $F\subseteq A^\ast$ be finite. The subset 
$$X_F=\{x\in A^\mathbb{Z} \mid \text{no word in}\; F\;\text{appears in}\ x \} $$
 is clearly a subshift. Such subshifts will be called shifts of finite type (SFT's in short).
\end{definition}

Thus SFT's are subshifts determined by the finitely many \textit{forbidden words} in $F$. 
When $w\in A^\ast$ belongs to $L(X_F)$ (equivalently $w$ has no subwords from $F$), we will say that $w$ is \textit{admissible}.

In some sense SFT's and minimal subshifts lie in the opposite spectrum. One aspect of this is the fact that if an SFT is minimal it is finite. Another aspect is that while  topological full groups of minimal subshifts are amenable (and hence cannot contain non-abelian free subgroups), topological full groups of SFT's can contain non-abelian free subgroups (see Section \ref{f2sft} Theorem \ref{thm:irreducible-sft-free-subgroup}).

\vspace{0.2cm}
If $X_F$ is an SFT with $\ell=\max\{|w|\mid w\in F\}$, define 
$$F'=\{w\in A^\ell \mid w\; \text{has a subword from}\; F\}.$$
It is clear that $X_F=X_{F'}$ and hence we may assume that an SFT is defined by a sef of forbidden words of uniform length $\ell$. Thus, a word $w\in A^\ast$ is admissible if and only if
all subwords of $w$ of length $\ell-1$ are admissible. In this case we will say that $X_F$
an $\ell-1$ \textit{step shift of finite type}. More precisely,

\begin{definition}
An SFT  is an $m$-step SFT if it is defined by a set of forbidden words of lengths $m+1$.
\end{definition}

As remarked above, infinite SFT's are never minimal. But they may be \textit{topologically transitive} meaning that there is at least one dense orbit. It can be observed that for subshifts topological transitivity is equivalent to the following:

\begin{definition}
A subshift $X\subseteq A^\mathbb{Z}$ is irreducible if for every $u,v\in L(X)$, there is
$w\in L(X)$ such that $uwv\in L(X)$.
\end{definition}

A uniform version of irreducibility is the notion of \textit{mixing}:
\begin{definition}
A subshift $X\subseteq A^\mathbb{Z}$ is topologically mixing if for every $u,v\in L(X)$ there is $N\ge 1$ such that for every $n\ge N$ there is $w\in A^n$ with $uwv\in L(X)$.
\end{definition}

\subsection{Approximating groups and finitely presented covers}\label{approx_seq} \mbox{}

\vspace{0.2cm}

Let $G$ be a finitely generated group given by a (not-necessarily finite) presentation
$$G=\left<x_1,\ldots,x_k \mid r_1,r_2,r_3,\ldots  \right>.$$
 Truncating  this presentation one obtains finitely presented groups 
$$G_n= \left< x_1,\ldots x_k \mid r_1,\ldots r_n\right>$$
which \textit{approximate} $G$. Let us call such a squence an \textit{approximating sequence}. In fact, $G_n$ converge to $G$ in the space of
$k$\textit{-generated groups}. (We will not give the definition and topology of the space of 
finitely generated groups since it will not be used in this article. Interested readers can
see \cite[Section 3]{BGdlH2013}).

A \textit{finitely presented cover of $G$} is a a finitely presented group $H$ together with en epimorphism $H\to G$. Such covers for a variety of groups were investigated in \cite{BGdlH2013}.

The following well known fact relates these notions:

\begin{lemma}\label{BGH}\cite[Proposition 3.3]{BGdlH2013}
Let $G$ and $G_n$ be as above. If $H$ is a finitely presented cover of $G$, then $H$ is a finitely presented cover of $G_n$ for large enough $n$.
\end{lemma}

A general method of showing that a finitely generated group $G $ is not finitely presented is the following: Suppose the group $G$ does not have a property $P$ and this property is inherited by homomorphic pre-images (i.e., whenever a group $H$ has a quotient with $P$ then $H$ has $P$, for example \textit{having a non-abelian free subgroup} is such a property).
Suppose and for some approximating sequence $G_n$ as above, the groups $G_n$ do not have property $P$. If one can find an approximating sequence $G_n$ with all $G_n$ having property $P$, then $G$ can't be finitely presented by the previous Lemma. Let us formulate this as a lemma when $P$ is the property of having non-abelian free subgroups:

\begin{lemma}
Let $G$ be a finitely generated group which does not have non-abelian free subgroups. If there is an approximating sequence $G_n$ such that each $G_n$ has non-abelian free subgroups, then every finitely presented cover of $G$ has non-abelian free subgroups. In particular, $G$ is not finitely  presented.
\end{lemma}


\section{Main Results and proofs}

\subsection{Free Subgroups and shifts of finite type} \label{f2sft} \mbox{}

\vspace{0.2cm}

In this section, we will show that topological full groups of infinite, irreducible SFT's contain non-abelian free subgroups. We first prove a criterion which ensures the existence of 
non-abelian free subgroups. The idea of this theorem was conceived during a private communication with H. Matui.

\begin{theorem}\label{thm:sft-free-product}
Let $(X,T)$ be an $m$-step SFT over the alphabet $A$. Suppose that there are  integers $N,n\ge m$, distinct words $u,v,w\in L_n(X)$ with the following property:

For every ordered pair of distinct words $a,b\in\{u,v,w\}$, there exists a word
	$\eta_{ab}\in A^{M-n}$,
	such that
	$a\eta_{ab}b\eta_{bc}c\in L(X)$
	for all
	$a,b,c\in\{u,v,w\}$
	with
	$a\neq b$
	and
	$b\neq c$.

	Then 
	$\fullgp{T}$
	contains a subgroup isomorphic to
	$\mathbb Z_2*\mathbb Z_2*\mathbb Z_2$.
	In particular, $\fullgp{T}$ contains a non-abelian free subgroup.
\end{theorem}

\begin{proof}
Define the clopen sets
	\begin{align*}
		U_u&=\left\{x\in X\mid x_{[0,n-1]}=u,\ x_{[M,M+n-1]}=v\ \text{or}\ x_{[M,M+n-1]}=w\right\}
		\\
		U_v&=\left\{x\in X\mid x_{[0,n-1]}=v,\ x_{[M,M+n-1]}=u\ \text{or}\ x_{[M,M+n-1]}=w\right\}
		\\
		U_w&=\left\{x\in X\mid x_{[0,n-1]}=w,\ x_{[M,M+n-1]}=u\ \text{or}\ x_{[M,M+n-1]}=v\right\}		
	\end{align*}

	By assumption, for each distinct $a,b\in\{u,v,w\}$, $a\eta_{ab}b\in L(X)$ and hence
	$[\bullet a\eta_{ab}b]\subseteq U_a$, so the set $U_a$ is nonempty. For each $a\in\{u,v,w\}$ define the homeomorphisms:
	$$
	\Delta_a(x)=
	\begin{cases}
		T^M(x), & x\in U_a,\\
		T^{-M}(x), & x\in T^M(U_a),\\
		x, & \text{otherwise}.
	\end{cases}
	$$
	 Clearly,
	$\Delta_u,\Delta_v,\Delta_w\in\fullgp{T}$ and
	$\Delta_u^2=\Delta_v^2=\Delta_w^2=1$.
	
	Consider the free product
	$$
	\Gamma=\mathbb Z_2*\mathbb Z_2*\mathbb Z_2
	=
	\langle u,v,w\mid u^2=v^2=w^2=1\rangle.
	$$
It follows that there is a homomorphism 
	$\mu\colon\Gamma\longrightarrow\fullgp{T}$
 with 
	$\mu(g)=\Delta_g$
	for all $g\in\{u,v,w\}$. We will show that $\mu$ is injective.

	Let
	$g_k g_{k-1}\cdots g_1$
	be a nonempty reduced word in $\Gamma$, where each
	$g_i\in\{u,v,w\}$
	and
	$g_{i+1}\neq g_i$
	for every $i=1,\ldots,k-1$. Choose
	$d\in\{u,v,w\}\setminus\{g_1,g_k\}$.
	By the hypothesis, the word
	$$
	g_1\eta_{g_1g_2}g_2\eta_{g_2g_3}g_3\cdots
	g_{k-1}\eta_{g_{k-1}g_k}g_k\eta_{g_kd}d
	$$
	is admissible in $X$. More precisely, since $X$ is an $m$-step shift of finite type, it is enough to check admissibility of the blocks of length $m+1$. Since
	$n\geq m$, these blocks occur inside words of the form
	$a\eta_{ab}b\eta_{bc}c,$
	which are admissible by assumption.
	Hence there exists  $z\in X$
	with
	$$
	z_{[0, n-1]}=g_1, \
	z_{[M, M+n-1]}=g_2, \
	\ldots, \
	z_{[(k-1)M, (k-1)M+n-1]}=g_k, \
	z_{[kM, kM+n-1]}=d.
	$$
	Now:
	\begin{align*}
		\mu(g_k\cdots g_1)(z)
		&=\Delta_{g_k}\cdots\Delta_{g_1}(z)\\
		&=\Delta_{g_k}\cdots\Delta_{g_2}\bigl(T^{M}(z)\bigr)\\
		&\ \vdots\\
		&=T^{kM}(z)\neq z
	\end{align*}
	as $z\in[\bullet g_1]$ and $T^{kM}(z)\in[\bullet d]$, while
	$[\bullet g_1]\cap[\bullet d]=\varnothing$. Therefore
	$\mu(g_k g_{k-1}\cdots g_1)\neq 1$
	and hence $\mu$ is injective.
	It follows that 
$
	\langle\Delta_u,\Delta_v,\Delta_w\rangle
	\cong
	\mathbb Z_2*\mathbb Z_2*\mathbb Z_2.
$
and thus	 
	$
	\left\langle
	(\Delta_u\Delta_v)^2,\,
	(\Delta_v\Delta_w)^2
	\right\rangle \cong F_2
	$. 
\end{proof}

In general an irreducible SFT does not have to be mixing. The next lemma show that every irreducible SFT can be decomposed into mixing pieces:

\begin{lemma}\label{lem:cyclic-mixing-decomposition}{\cite[Theorem~2.1.13(2)]{Putnam2014}}, {\cite[Propositions~4.5.6 and~4.5.10]{LM1995}}
	
	Let $(X,T)$ be an irreducible $N$-step SFT.
	Then there exist an integer $d\geq 1$ and a non-empty clopen partition $	X=C_0\sqcup C_1\sqcup\cdots\sqcup C_{d-1}$ such that
	 $$	T(C_i)=C_{i+1}\quad\text{for }0\leq i<d-1,
	\quad
	T(C_{d-1})=C_0,
	$$
	and, for every $i=0,\ldots,d-1$, the restriction
	$T^d\big|_{C_i}\colon C_i\longrightarrow C_i$
	is topologically mixing. Equivalently, for every pair of nonempty open
	sets $U,V\subseteq C_i$, there exists $N\geq 1$ such that
	$T^{nd}(U)\cap V\neq\varnothing$
	for every $n\geq N$.
\end{lemma}

The next lemma show that the hypotheses of Theorem \ref{thm:sft-free-product} are fulfilled for irreducible SFT's.

\begin{lemma}
	\label{lem:three-cylinders-in-cyclic-component}
	Let $(X,T)$ be an infinite irreducible $N$-step SFT over the alphabet $A$ and let 
	$X=C_0\sqcup C_1\sqcup\cdots\sqcup C_{d-1}$
	be the cyclic decomposition given by Lemma~\ref{lem:cyclic-mixing-decomposition}. Fix $i\in\{0,\ldots,d-1\}$. Then there exist an integer $n\geq N$, three distinct words
	$u,v,w\in L_n(X)$
	with
	$$
	[\bullet u],[\bullet v],[\bullet w]\subseteq C_i,
	$$
	and an integer $M\geq n$, divisible by $d$, such that for every ordered pair of distinct words $a,b\in\{u,v,w\}$, there exists a word
	$\eta_{ab}\in A^{M-n}$
	satisfying
	$$
	a\eta_{ab}b\eta_{bc}c\in L(X)
	$$
	for all $a,b,c\in\{u,v,w\}$ with $a\neq b$ and $b\neq c$.
\end{lemma}

\begin{proof}
	By Lemma~\ref{lem:cyclic-mixing-decomposition}, the restriction
	$T^d|_{C_i}$ is topologically mixing. Since $T$ cyclically permutes the
	sets $C_0,\ldots,C_{d-1}$ and $X$ is infinite, the set $C_i$ is
	infinite. We may therefore choose an integer $\ell\geq N$ and two distinct
	words $\alpha,\beta\in L_{\ell}(X)$ such that
	$
	[\bullet\alpha],[\bullet\beta]\subseteq C_i.
	$ For every ordered pair $a,b\in\{\alpha,\beta\}$, the topological mixing of
	$T^d|_{C_i}$ gives an integer $k_{ab}\geq 1$ such that
	$$
	T^{kd}([\bullet a])\cap[\bullet b]\neq\varnothing
	$$
	for every $k\geq k_{ab}$. Choose $k$ larger than the four integers
	$k_{ab}$ and sufficiently large so  that $kd\geq\ell$, and put $L=kd$.
	Then, for every $a,b\in\{\alpha,\beta\}$, there exists a word
	$\theta_{ab}\in A^{L-\ell}$ such that
	$$
	a\theta_{ab}b\in L(X).
	$$	Now choose a sequence
	$
	a_0a_1a_2\cdots\in\{\alpha,\beta\}^{\mathbb{N}}
	$
	which is not eventually periodic. Since $X$ is an $N$-step shift of
	finite type and $\ell\geq N$, the admissible words
	$$
	a_j\theta_{a_ja_{j+1}}a_{j+1}
	$$
	can be successively joined along their common occurrences of $a_{j+1}$.
	Thus the right-infinite word
	$$
	a_0\theta_{a_0a_1}a_1\theta_{a_1a_2}a_2
	\theta_{a_2a_3}\cdots
	$$
	is admissible and has a left extension to a point $\xi\in X$. Moreover,
	$\xi\in[\bullet a_0]\subseteq C_i$, and the word $a_j$ begins at
	coordinate $jL$ for every $j\geq 0$.
	
	The right-infinite sequence determined by $\xi$ is not eventually
	periodic. Indeed, if it were eventually periodic with period $p$, then,
	for all sufficiently large $j$, the blocks of length $\ell$ beginning
	at the coordinates $jL$ would be periodic in $j$, with period
	$p/\gcd(p,L)$. Hence the sequence $(a_j)_{j\geq 0}$ would be eventually
	periodic, contrary to its choice.
	
	The right-infinite sequences determined by
	$\xi, \ T^d(\xi), \ T^{2d}(\xi)$
	are pairwise distinct. Indeed, equality of two of these right-infinite sequences would make the right-infinite sequence determined by $\xi$ eventually periodic. We may therefore choose $n\geq N$ such that the words
	$$
	u=\xi_{[0,n-1]},\qquad v=\xi_{[d,d+n-1]},\qquad w=\xi_{[2d,2d+n-1]}
	$$
	are pairwise distinct.
	
	Since $\xi,T^d(\xi),T^{2d}(\xi)\in C_i$ and $n\geq N$, we have
	$[\bullet u],[\bullet v],[\bullet w]\subseteq C_i.$
	In particular, these are pairwise disjoint nonempty cylinder sets.
	
	By Lemma~\ref{lem:cyclic-mixing-decomposition}, the restriction $T^d|_{C_i}$ is topologically mixing. Hence, for every ordered pair of distinct words $a,b\in\{u,v,w\}$, there exists $k_{ab}\geq 1$ such that
	$$
	T^{kd}([\bullet a])\cap[\bullet b]\neq\varnothing
	$$
	for every $k\geq k_{ab}$. Choose $k$ larger than all six integers $k_{ab}$ and sufficiently large so that $kd\geq n$, and put
	$M=kd.$
	Then $M\geq n$, $d\mid M$, and
	$$
	[\bullet a]\cap T^{-M}([\bullet b])\neq\varnothing
	$$
	for every ordered pair of distinct words $a,b\in\{u,v,w\}$.
	For each such pair, choose
	$$
	x^{ab}\in[\bullet a]\cap T^{-M}([\bullet b]).
	$$
	By the definition of these cylinders,
	$x^{ab}_{[0,n-1]}=a
	\ \text{and} \
	x^{ab}_{[M,M+n-1]}=b.$
	Since $M\geq n$, we may define
	$$
	\eta_{ab}=x^{ab}_{[n,M-1]}.
	$$
	Then $\eta_{ab}\in A^{M-n}$, and the block of $x^{ab}$ on the
	coordinates $0,\ldots,M+n-1$ is
	$a\eta_{ab}b.$
	Therefore
	$a\eta_{ab}b\in L(X).$
	
	Now let $a,b,c\in\{u,v,w\}$ satisfy $a\neq b$ and $b\neq c$.
	The element $x^{ab}$ contains the word $a\eta_{ab}b$, with the
	occurrence of $b$ beginning at coordinate $M$. On the other hand,
	  $x^{bc}$ contains the word $b\eta_{bc}c$, with the
	occurrence of $b$ beginning at coordinate $0$. Put
	$$
	y=T^{-M}(x^{bc}).
	$$
	Then the occurrence of $b$ in $y$ begins at coordinate $M$, while
	the occurrence of $c$ begins at coordinate $2M$. Hence
	$$
	x^{ab}_{[M,M+n-1]}=y_{[M,M+n-1]}=b.
	$$
	Thus $x^{ab}$ and $y$ agree on the $n$ consecutive coordinates
	$M,\ldots,M+n-1$. Since $n\geq N$, they agree there on an
	$N$-block.
	
	We may therefore join the left-hand part of $x^{ab}$ to the
	right-hand part of $y$ along this common $N$-block. More precisely,
	define $z\in A^{\mathbb Z}$ by
	$$
	z_j=
	\begin{cases}
		x^{ab}_j, & j\leq M+n-1,\\
		y_j, & j\geq M+n.
	\end{cases}
	$$

	Every block of $z$ of length $N+1$ either occurs in $x^{ab}$, occurs
	in $y$, or crosses the joining position. In the last case, it agrees
	with a block of $y$, because $x^{ab}$ and $y$ agree on the
	preceding $N$ coordinates. Since $X$ is an $N$-step shift of
	finite type, it follows that $z\in X$.
	
	The point $z$ contains $a$ at coordinate $0$, $b$ at coordinate
	$M$, and $c$ at coordinate $2M$. More precisely,
	$z_{[0,2M+n-1]}=a\eta_{ab}b\eta_{bc}c.$
	Therefore
	$
	a\eta_{ab}b\eta_{bc}c\in L(X)
	$
	as required.
\end{proof}

We are now ready to state and prove the main theorem of this section:
\begin{theorem}
	\label{thm:irreducible-sft-free-subgroup}
	Let $(X,T)$ be an infinite irreducible SFT. Then
	$\fullgp{T}$ contains a subgroup isomorphic to
	$\mathbb Z_2*\mathbb Z_2*\mathbb Z_2$. In particular,
	$\fullgp{T}'$ contains a non-abelian free subgroup.
\end{theorem}

\begin{proof}
	Let $X$ be an $N$-step shift of finite type.
	By Lemma \ref{lem:three-cylinders-in-cyclic-component}, there exist three
	distinct words $u,v,w\in L_n(X)$, where $n\geq N$, and an integer
	$M\geq n$ satisfying the hypotheses of
	Theorem~\ref{thm:sft-free-product}. Therefore there exist involutions
	$\Delta_u,\Delta_v,\Delta_w\in\fullgp{T}$ such that
$
	\langle\Delta_u,\Delta_v,\Delta_w\rangle
	\cong
	\mathbb Z_2*\mathbb Z_2*\mathbb Z_2.
	$
	and hence
	$
	F_2\cong \left\langle
	(\Delta_u\Delta_v)^2,\,
	(\Delta_v\Delta_w)^2
	\right\rangle
	\le \fullgp{T}'.
	$
	
\end{proof}

\subsection{Approximating minimal subshifts by SFT's} \label{approxsft} \mbox{}

\vspace{0.2cm}
Let $(X,T)$ be a subshift over an alphabet $A$. For $N\ge 1$, let $A_{\le N}=\bigcup_{i=1}^N A^i$ and $L_{\le N}(X)=\bigcup_{i=1}^N L_i(X)$. For $N\ge 1$, let $F_N=A_{\le N+1} \setminus L_{\le N+1}$ and consider the SFT
$$X_N=X_{F_N}=\{ x\in A^\mathbb{Z} \mid \text{no word in}\; F_N\;\text{appears in}\ x  \}. $$
In  words, $X_N$ is an SFT for which the forbidden words are exactly the words of length at most $N+1$ that do not appear in sequences of $X$. Clearly, $X\subseteq X_N$ and $X_{N+1}\subseteq X_N$ for every $N\ge 1$.  Also note that $X_N$ is an $N$-step SFT with 
$$L_N(X_N)=L_N(X)\; \text{and} \; L_{\le N}(X_N)=L_{\le N}(X). $$

The next lemma shows that for a transitive subshift (in particular for a minimal one), these approximating SFT's are irreducible.

\begin{lemma}\label{lem:minimal-approximants-irreducible}
	Let $X\subseteq A^{\mathbb Z}$ be a transitive subshift.
	Then, for every $N\geq 1$, the SFT $X_N$
	defined above is irreducible.
\end{lemma}

\begin{proof}

	Let $u,v\in L(X_N)$. We shall construct a word joining $u$ to
	$v$ in $X_N$.
	
	Since $u\in L(X_N)$, we may extend $u$ to the right by an
	$N$-block $w\in L_N(X_N)$ such that
	$uw\in L(X_N).$
	Similarly, we may extend $v$ to the left by an $N$-block
	$t\in L_N(X_N)$ such that
	$tv\in L(X_N).$
	More explicitly, choose $x,y\in X_N$ satisfying
	$$
	x_{[0,|u|-1]}=u
	\qquad\text{and}\qquad
	y_{[0,|v|-1]}=v,
	$$
	and put
	$$
	w=x_{[|u|,|u|+N-1]},
	\quad
	t=y_{[-N,-1]}.
	$$
	Since
	$w,t\in L_N(X_N)=L_N(X)$
	and $X$ is transitive, there exists a word $\eta$ such that
	$$
	w\eta t\in L(X).
	$$
	Since $X\subseteq X_N$, it follows that
	$$
	w\eta t\in L(X_N).
	$$	We now have three admissible words
	$uw, \ w\eta t, \ tv$
	in $L(X_N)$. The first two words agree on the $N$-block $w$, and
	the last two agree on the $N$-block $t$. Since $X_N$ is an
	$N$-step SFT, they may be joined along these common
	$N$-blocks. Thus
	$$
	uw\eta tv\in L(X_N).
	$$
	Therefore the word
	$w\eta t$
	connects $u$ to $v$ in $X_N$. Since $u,v\in L(X_N)$ were
	arbitrary, $X_N$ is irreducible.
\end{proof}

\subsection{Approximating topological full groups} \mbox{}

\vspace{0.2cm}

Let $(X,T)$ be an infinite minimal subshift over an alphabet $A$. As mentioned in the preliminaries (Lemma \ref{norepeat}) we can assume without loss of generality that the subshift satisfies the conclusion of Lemma \ref{norepeat}. Throughout this section we have this convention and we skip mentioning it repeatedly.

Let $\left<x_{(w,1)}, w\in L_3(X) \mid R \right>$ be the presentation of $\fullgp{T}'$ given by Grigorchuk and Medynets (Theorem \ref{gmpres}). Recall that the relations $R$ are parametrized by cylinder sets of the form $(w,i)$ where $w\in L(X)$ and $1\le i \le |w|-1$ and also cylinder partitions of $(w,i)$ into smaller cylinder sets of this form.

For every $N\geq 5$, let $R_N$ denote the finite collection of all instances of the relations in $R$  which involve only admissible cylinder words of length at most $N+1$. We define the approximating sequence (as in Section \ref{approx_seq}) 
$$G_N=\left<x_{(w,1)}, w\in L_3(X) \mid R_N \right>. $$

The aim of this section is to show that  $G_N$ has non-abelian free subgroups, using ideas from Sections \ref{f2sft} and \ref{approxsft}.

\vspace{0.2cm}

For each $N\ge 1$, recall the approximating SFT's $X_N$ of Section \ref{approxsft}. $X_N$ is the SFT whose forbidden words are the words of length at most $N+1$ which do not appear in sequences in $X$.

 To distinguish the shift map $T$ on $X$ and $X_N$, we will write $T_N:X_N \to X_N$ for the shift map on $X_N$. Thus, the minimal subshift $(X,T)$ is approximated by the  SFT's $(X_N,T_N)$, which are irreducible by Lemma \ref{lem:minimal-approximants-irreducible}.

 In the remainder we assume that  $N\ge 5$. By construction,
$L_{\leq N+1}(X_N)=L_{\leq N+1}(X).$ and  since $N\geq 5$,
$L_{\leq 3}(X_N)=L_{\leq 3}(X).$

\vspace{0.2cm}

For
$
w\in L_3(X)=L_3(X_N),
$
define the local 3 cycle
$
\sigma_N(w,1)\in \fullgp{T_N}
$
by
$$
\sigma_N(w,1)(z)=
\begin{cases}
	T_N(z), & z\in (w,1)\sqcup T_N((w,1)),\\
	T_N^{-2}(z), & z\in T_N^2((w,1)),\\
	z, & \text{otherwise}.
\end{cases}
$$
As before, 
$(w,1)=T_N[\bullet w]\subseteq X_N.$
Since $X$ satisfies Lemma \ref{norepeat} and
$L_{\leq N+1}(X_N)=L_{\leq N+1}(X)$,
the cylinders needed to define these local 3 cycles are pairwise disjoint in
$X_N$. Hence $\sigma_N(w,1)$ is well-defined.

 Let us call cylinder sets $(w,i),(v,j)$ \textit{3 disjoint}, whenever  $(w,i),(w,i+1),w(i+2),(v,j),(v,j+1),(v,j+2)$ are disjoint.

\begin{theorem}
	\label{thm:GN-to-TN}
	 There is a well defined homomorphism $\Phi_N\colon G_N\longrightarrow \fullgp{T_N}$ which sends
	$x(w,1)\mapsto \sigma_N(w,1)$ for every $w\in L_3(X)$.

\end{theorem}
\begin{proof}
	Define $\Phi_N$ on the free group generated by
	$
	\{x(w,1)\mid w\in L_3(Y)\}
	$
	by
	$$\Phi_N(x(w,1))=\sigma_N(w,1).$$
	We will show that every relation in $R_N$ is sent to the identity.
	Since
	$L_{\leq N+1}(X_N)=L_{\leq N+1}(X),$
	the cylinder intersections, cylinder partitions, and $3$-disjointness
	conditions which appear in the relations $R_N$ are the same in $X$ and
	in $X_N$. Therefore the recursive definitions of the formal symbols
	$x(w,i)$ agree with the corresponding local three-cycles
	$\sigma_N(w,i).$
	In other words, for every formal expression $x(w,i)$ appearing in $R_N$,
	the image under $\Phi_N$ is
	$\Phi_N(x(w,i))=\sigma_N(w,i).$
	
	It remains to check the defining relations. The element $\sigma_N(v,i)$ is a
	three-cycle on
	$(v,i), \ (v,i+1), \ (v,i+2),$
	and hence
	$\sigma_N(v,i)^3=1.$
	Also,
	$$
	\bigl(\sigma_N(v,i)\sigma_N(v,i+1)\bigr)^2=1
	$$
	and
	$$
	\sigma_N(v,i+1)=\sigma_N(v,i)*\sigma_N(v,i+2)
	$$
	are the standard relations among these local three-cycles.
	
	If $(v,i)$ and $(w,j)$ are $3$-disjoint, then the supports of
	$\sigma_N(v,i)$ and $\sigma_N(w,j)$
	are disjoint. Hence
	$$
	[\sigma_N(v,i),\sigma_N(w,j)]=1.
	$$
	Finally, if $C$ is a  cylinder partition of
	$(v,i)$, then
	$$
	\sigma_N(v,i)
	=
	\prod_{(s,k)\in C}
	\sigma_N(s,k),
	$$
	because the atoms of the partition are pairwise disjoint and the local
	three-cycle acts independently on each atom.
	
	Thus every defining relation in $R_N$ is satisfied by the corresponding
	elements in $\fullgp{T_N}$. Thus, the assignment
	$x(w,1)\longmapsto \sigma_N(w,1)$
	extends to a well-defined homomorphism
	$\Phi_N\colon G_N\longrightarrow \fullgp{T_N}.$
\end{proof}

Let
$H_N=\Phi_N(G_N)\leq \fullgp{T_N}$. Since $X\subseteq X_N$ and $X$ is infinite, the shift $X_N$ is also infinite. Thus by Lemma \ref{lem:minimal-approximants-irreducible} and Theorem \ref{thm:irreducible-sft-free-subgroup}, $\fullgp{T_N}$ contains a non-abelian free group. The   next lemma shows that this free subgroup   belongs to $H_N$.

\begin{lemma}
	\label{lem:squares-in-HN}
	Let $u,v,w\in L_n(X_N)$ and $M\geq n$ satisfy the hypotheses of
	Theorem \ref{thm:sft-free-product}, and let
	$\Delta_u,\Delta_v,\Delta_w\in\fullgp{T_N}$
	be the corresponding involutions constructed in its proof. Then
	$
	(\Delta_u\Delta_v)^2, (\Delta_v\Delta_w)^2
	\in H_N$. Thus, for $N\ge 5$, $G_N$ contains a non-abelian free subgroup.
\end{lemma}

\begin{proof}
	For each $x,y\in\{u,v,w\}$, since the elements $\Delta_x \Delta_y$ are involutions, we have 
	$(\Delta_u\Delta_v)^2=[\Delta_u,\Delta_v].$
	Choose a finite clopen partition which simultaneously refines the supports
	of $\Delta_u$ and $\Delta_v$ and their images. On the resulting clopen
	atoms, both $\Delta_u$ and $\Delta_v$ act as products of disjoint
	two-cycles. Hence $(\Delta_u\Delta_v)^2$ acts as an even permutation of
	these atoms. Such a permutation is a product of three-cycles. Explicitly,
	for clopen atoms $A,B,C,D$, the elementary identities
	$$
	(A\,B)(B\,C)=(A\,B\,C) \quad \text{and} \quad
	(A\,B)(C\,D)=(A\,C\,B)(A\,C\,D)
	$$
	show that every product of an even number of two-cycles can be written as
	a product of three-cycles.
	
	We now verify that the three-cycles obtained above belong to $H_N$.
	Refine the clopen atoms further so that, on each atom, every partial
	product occurring in the word
	$$
	[\Delta_u,\Delta_v]
	=\Delta_u\Delta_v\Delta_u\Delta_v
	$$
	acts as a fixed power of $T_N$. Since all these maps preserve
	$T_N$-orbits, the resulting finite permutations may be considered
	separately on finite collections of atoms contained in the same
	$T_N$-orbit. After adjoining the finitely many intermediate translates
	when necessary, order these atoms according to their positions in the
	orbit. The alternating group on such a finite ordered collection is
	generated by the adjacent three-cycles
	$
	(j\;\;j+1\;\;j+2).
	$
	Each of these adjacent three-cycles is a local three-cycle associated
	to a clopen set $C$, with support
	$
	C,T_N(C),T_N^2(C)
	$. Refining $C$ into cylinder sets and using the cylinder-refinement
	identity, this local three-cycle becomes a product of
	cylinder three-cycles. Moreover, since $N\geq 5$ and
	$$
	L_{\leq N+1}(X_N)
	=L_{\leq N+1}(X),
	$$
	the single-letter cylinders in $X_N$ are disjoint from their first
	four translates, exactly as in $X$. Hence the shortening and
	dot-moving identities established in Section \ref{symbdynand}  apply verbatim to these
	cylinder three-cycles; these identities are local and do not require
	minimality. Consequently, every such cylinder three-cycle is a word in
	$$
 	\sigma_N(q,1), q\in L_3(X_N)
	=L_3(X).
	$$
	Therefore every three-cycle appearing in the decomposition of
	$(\Delta_u\Delta_v)^2$ belongs to $H_N$, and hence
	$
	(\Delta_u\Delta_v)^2\in H_N.
	$
	The same argument, with $u,v$ replaced by $v,w$, gives
	$
	(\Delta_v\Delta_w)^2\in H_N
	$.
	
	It follows that the image of $\Phi_N\colon G_N\longrightarrow \fullgp{T_N}$ contains a non-abelian free subgroup. Thus $G_N$ contains a non-abelian free subgroup.
\end{proof}

We are now ready to prove the main theorem.

\begin{theorem}\label{mainthm}
Let $(X,T)$ be an infinite minimal subshift. If $G$ is a finitely presented group
which maps onto $\fullgp{T}'$, then $G$ contains a non-abelian free subgroup. In 
particular, $\fullgp{T}'$ is not finitely presented.
\end{theorem}

\begin{proof}
Without loss of generality we can assume that $(X,T)$ satisfies Lemma \ref{norepeat}. By Lemma
\ref{BGH}, $G$ maps onto $G_N$ for large enough $N$. By Theorem \ref{thm:GN-to-TN}, $G_N$ contains a non-abelian free subgroup, which implies that $G$ contains a non-abelian free subgroup. If $\fullgp{T}'$ was finitely presented, it would contain a non-abelian free subgroup which is not possible since $\fullgp{T}'$ is amenable by Theorem \ref{ref} item (5).
\end{proof}

\bibliography{fullgpspaper}{}
\bibliographystyle{alpha}

\end{document}